\documentclass{gtart}

\input gtoutput
\volumenumber{4}\papernumber{4}\volumeyear{2000}
\pagenumbers{149}{170}\published{11 March 2000}
\proposed{Ralph Cohen}\seconded{Gunnar Carlsson, Yasha Eliashberg}
\received{30 July 1999}
\accepted{29 February 2000}

\usepackage{amssymb}

\newtheorem{theorem}{Theorem}[section]
\newtheorem{proposition}[theorem]{Proposition}
\newtheorem{lemma}[theorem]{Lemma}
\newtheorem{corollary}[theorem]{Corollary}
\theoremstyle{definition}

\newcommand{\R}{{\mathbb R}}
\newcommand{\Z}{{\mathbb Z}}
\newcommand{\C}{{\mathbb C}}

\renewcommand{\Im}{{\mbox{Imm}}}
\newcommand{\Sq}{{\mbox{Sq}}}

\begin{document}

\title{Double point self-intersection surfaces of immersions}
%\shorttitle{Double space surfaces}
\covertitle{Double point self-intersection\\surfaces of immersions}
\author{Mohammad A Asadi-Golmankhaneh\\Peter J Eccles}
\shortauthors{Asadi-Golmankhaneh and Eccles}
\address{Department of Mathematics, University of Urmia\\
PO Box 165, Urmia, Iran}
\secondaddress{Department of Mathematics, University of
Manchester\\Manchester, M13 9PL, UK}
\email{pjeccles@man.ac.uk}
\primaryclass{57R42}
\secondaryclass{55R40, 55Q25, 57R75}
\keywords{immersion, Hurewicz homomorphism, spherical
class, Hopf invariant, Stiefel--Whitney number}
\asciikeywords{immersion, Hurewicz homomorphism, spherical
class, Hopf invariant, Stiefel-Whitney number}

\begin{abstract}
A self-transverse immersion of a smooth manifold $M^{k+2}$
in $\R^{2k+2}$ has a double point self-intersection
set which is the
image of an immersion of a smooth surface, the double point
self-intersection surface.  We prove that
this surface may have odd Euler characteristic if and only
if $k\equiv 1\bmod 4$ or $k+1$ is a power of 2.  This
corrects a previously published result by Andr\'{a}s
Sz\H{u}cs \cite{szu93}.

The method of proof is to evaluate the Stiefel--Whitney
numbers of the double point self-intersection surface.  By
the methods of \cite{asa99} these numbers can be read off
from the Hurewicz image $h(\alpha)\in
H_{2k+2}\Omega^{\infty}\Sigma^{\infty}MO(k)$ of the element
$\alpha\in\pi_{2k+2}\Omega^{\infty}\Sigma^{\infty}MO(k)$
corresponding to the immersion under the Pontrjagin--Thom
construction.
\end{abstract}

\asciiabstract{A self-transverse immersion of a smooth manifold 
M^{k+2} in R^{2k+2} has a double point self-intersection set which is
the image of an immersion of a smooth surface, the double point
self-intersection surface.  We prove that this surface may have odd
Euler characteristic if and only if k is congruent to 1 modulo 4 or
k+1 is a power of 2.  This corrects a previously published result by
Andras Szucs.

The method of proof is to evaluate the Stiefel-Whitney numbers of the
double point self-intersection surface.  By earier work of the authors
these numbers can be read off from the Hurewicz image h(\alpha) in
H_{2k+2}\Omega^{\infty}\Sigma^{\infty}MO(k) of the element \alpha in
\pi_{2k+2}\Omega^{\infty}\Sigma^{\infty}MO(k) corresponding to the
immersion under the Pontrjagin-Thom construction.}

\maketitlepage

\section{Introduction}

Let $f\co M^{n-k} \looparrowright \R^n$ be a
self-transverse immersion of a compact closed smooth
($n-k$)--dimensional manifold in $n$--dimensional Euclidean
space ($0 < k \leq n$).  A point of $\R^n$ is an {\em
$r$--fold self-intersection point} of the immersion if it is the
image under $f$ of (at least) $r$ distinct points of the
manifold.  The self-transversality of $f$ implies that the
set of $r$--fold self-intersection points (the {\em $r$--fold
self-intersection set})
is itself the image of an immersion
\[\theta_r(f)\co \Delta_r(f)\looparrowright \R^n\]
of a compact manifold
$\Delta_r(f)$, the {\em $r$--fold self-intersection manifold}, of
dimension $n-rk$, although in general this immersion will not
be self-transverse.  It is natural to ask what can be said
about this $r$--fold intersection manifold: which manifolds
can arise for each value of $n$, $k$ and $r$?

The simplest case is when $n=rk$ so that the self-intersection
manifold is a finite set of points and this case was the
first to be considered in detail (see \cite{ban74},
\cite{ecc80}, \cite{ecc81}, \cite{ecc93}, \cite{lan84}).

More recently cases of
higher dimensional self-intersection manifolds have been
considered.
Andr\'{a}s Sz\H{u}cs was one of the first to do so and in
\cite{szu93} considered the simplest case of interest, when
the double point self-intersection manifold is a surface; this arises for
$n=2k+2$.  
In this paper we return to this case using
different methods to Sz\H{u}cs.  Our main result is the
following.

\begin{theorem}\label{1.1} For $k\geq 1$, there exists an
immersion $f\co M^{k+2}\looparrowright \R^{2k+2}$ with
double point self-intersection
manifold of odd Euler characteristic if and only
if $k\equiv 1 \bmod4$ or $k+1$ is a power of $2$.
\end{theorem}

This result should be contrasted with Sz\H{u}cs' result
which asserted that double point self-intersection manifolds of odd Euler
characteristic can occur only if $k\equiv 1\bmod 4$.
Sz\H{u}cs's approach used differential topology and the
argument in the case $k\equiv 3\bmod 4$ was particularly
delicate.

Our approach uses algebraic topology and in particular the
correspondence between bordism groups and homotopy groups of
Thom complexes. In \cite{asa99} we described a general
approach to these problems which gives a method for
determining the
bordism class of the self-intersection manifolds of any
immersion: the unoriented bordism class of a manifold is determined by
its Stiefel--Whitney numbers and the Stiefel--Whitney numbers of
the self-intersection manifolds of an immersion can be read
off from certain homological information about the
immersion.  For double point self-intersection surfaces
the situation is particularly
simple since there are only two bordism classes:
a compact surface is a boundary if its Euler
characteristic is even and is is a non-boundary bordant to the projective
plane if its Euler characteristic is odd.

The paper is organized as follows.  In Section~2 we recall
the results needed from \cite{asa99}, establish our basic
notation and outline the proof of the theorem covering the
steps which apply in all cases.  In Section~3 we complete the
proof in the easiest case of $k$
even and this is followed by the proofs for $k\equiv 1\bmod 4$
in Section~4 and
$k\equiv 3 \bmod 4$ in Section~5. Finally, in Section~6, we
comment on the relationship between our results and those in
Sz\H{u}cs' paper \cite{szu93}.  Almost everything in that paper is
confirmed by our methods.

\rk{Acknowledgements} Most of the results in this paper are
contained in the first author's thesis \cite{asa98} which also contains
other applications of these methods.  He was supported by the
University of Urmia and the Ministry of Culture and Higher
Education of the Islamic Replublic of Iran during his time
as a student at the University of Manchester.

\section{The Stiefel Whitney characteristic numbers of the double point
self-intersection manifold}

Let $\Im(n-k,k)$ denote the group of bordism classes of
immersions $M^{n-k}\looparrowright \R^n$ of compact closed
smooth manifolds in Euclidean $n$--space.
Details of (co)bordism in this setting have been given by
R~Wells in \cite{wel66}.  By general position every
immersion is regularly homotopic (and so bordant) to a
self-transverse immersion and so each element of
$\Im(n-k,k)$ can be represented by a self-transverse
immersion.  In the same way bordisms between self-transverse
immersions can be taken to be
self-transverse; it is clear that such a bordism will
induce a bordism of the immersions of the double point self-intersection
manifolds so that $f\mapsto \theta_2(f)$ induces a double
point self-intersection map
\[ \theta_2\co \Im(n-k,k) \rightarrow \Im(n-2k,2k).\]

Wells shows how $\Im(n-k,k)$ may be described as a stable
homotopy group.
Let $MO(k)$ denote the Thom complex of a universal
$O(k)$--bundle $\gamma^k\co EO(k)\rightarrow BO(k)$
(see \cite{mil74} for basic material on vector bundles,
Thom complexes and bordism theory).
Using the Pontrjagin--Thom construction, Wells describes an isomorphism
\[\phi\co\Im(n-k,k)\cong\pi_n^SMO(k).\]

Write $QX$ for the direct limit
$\Omega^{\infty}\Sigma^{\infty}X = \lim \Omega^n \Sigma^n
X$,
where $\Sigma$ denotes the reduced suspension functor and
$\Omega$ denotes the loop space functor.  By the
adjointness of the functors $\Sigma$ and $\Omega$,
$\pi^S_nMO(k)\cong \pi_nQMO(k)$.
We consider the
$\Z/2$--homology Hurewicz homomorphism \[h\co
\pi_n^SMO(k)\cong \pi_nQMO(k) \rightarrow H_nQMO(k) =
H_n(QMO(k);\Z/2).\]  (Throughout this paper we use $H_*$ and
$H^*$ to denote homology and cohomology with $\Z/2$
coefficients.)

The main result of \cite{asa99} describes how, for a
self-transverse
immersion $f\co M^{n-k}\looparrowright\R^n$ corresponding
to $\alpha\in\pi_n^SMO(k)$, the Hurewicz image $h(\alpha)\in
H_nQMO(k)$
determines the normal Stiefel--Whitney numbers of the
self-intersection manifolds $\Delta_r(f)$.

To state this result in the case of double point self-intersection
manifolds we need some
preliminaries.

The {\em quadratic construction} on a
pointed space $X$ is defined to be
\[D_2X=X\wedge X\rtimes_{\Z/2} S^{\infty}
= X\wedge X\times_{\Z/2} S^{\infty}/\{*\}\times_{\Z/2}
S^{\infty},\]
where the non-trivial element of the group $\Z/2$ acts on
$X\wedge X$ by permuting the co-ordinates and on the
infinite sphere $S^{\infty}$ by the antipodal action.  There
is a natural map
\[h^2\co QX \rightarrow QD_2X\]
known as a {\em stable James--Hopf map} which induces a
{\em stable Hopf invariant}
\[h^2_*\co \pi^S_nX \rightarrow \pi^S_nD_2X\]
(see \cite{bar74} and \cite{sna74}).

If the self-transverse immersion $f\co
M^{n-k}\looparrowright \R^n$ corresponds to the element
$\alpha\in\pi^S_nMO(k)$, then the immersion of the double point
self-intersection manifold
$\theta_2(f)\co \Delta_2(f)\looparrowright\R^n$
corresponds to the element $h^2_*(\alpha)\in \pi_n^S
D_2MO(k)$ given by the
stable Hopf invariant (see \cite{kos78}, \cite{sz76a},
\cite{sz76b}, \cite{vog74}).

The immersion $\theta_2(f)$ corresponds to an
element in the stable homotopy of $D_2MO(k)$ because the
immersion of the double point self-intersection manifold automatically
acquires additional structure on its normal bundle, namely
that at each point $f(x_1)=f(x_2)$ the normal
$2k$--dimensional space is decomposed as the direct sum of
the two (unordered) $k$--dimensional normal spaces of $f$ at
the points $x_1$ and $x_2$.  The universal bundle for this
structure is
\[\gamma^k\times\gamma^k\times_{\Z/2}1\co
EO(k)\times EO(k)\times_{\Z/2}S^{\infty}
\rightarrow BO(k)\times BO(k)\times_{\Z/2}S^{\infty}\]
which has the Thom complex $D_2MO(k)$.

Forgetting this additional structure on the immersion
corresponds to applying the
map
\[\xi_*\co \pi^S_nD_2MO(k)\rightarrow \pi^S_nMO(2k)\]
induced by the map of Thom complexes $\xi\co
D_2MO(k)\rightarrow MO(2k)$ which comes from the map $BO(k)\times
BO(k) \times_{\Z/2} S^{\infty} \rightarrow BO(2k)$ classifying
the bundle $\gamma^k\times\gamma^k\times_{\Z/2}1$.  Thus, we
have the following commutative diagram.

\begin{equation} \setlength{\unitlength}{3cm}
\begin{picture}(4,0.35)(0,0.35)
\put(1,0){\makebox(0,0){$\pi^S_nMO(k)$}}
\put(2,0){\makebox(0,0){$\pi^S_nD_2MO(k)$}}
\put(3,0){\makebox(0,0){$\pi^S_nMO(2k)$}}
\put(1,0.7){\makebox(0,0){$I(n-k,k)$}}
\put(3,0.7){\makebox(0,0){$I(n-2k,2k)$}}
\put(1.4,0.7){\vector(1,0){1.2}}
\put(1.35,0){\vector(1,0){0.2}}
\put(2.4,0){\vector(1,0){0.2}}
\put(1,0.55){\vector(0,-1){0.4}}
\put(3,0.55){\vector(0,-1){0.4}}
\put(2,0.75){\makebox(0,0)[b]{$\theta_2$}}
\put(1.05,0.35){\makebox(0,0)[l]{$\phi$}}
\put(0.95,0.35){\makebox(0,0)[r]{$\cong$}}
\put(2.95,0.35){\makebox(0,0)[r]{$\cong$}}
\put(3.05,0.35){\makebox(0,0)[l]{$\phi$}}
\put(1.45,0.05){\makebox(0,0)[b]{$h^2_*$}}
\put(2.5,0.05){\makebox(0,0)[b]{$\xi_*$}}
\end{picture}
\end{equation}
\vspace{1cm}

Turning now to homology, observe that, by adjointness, the stable James--Hopf
map $h^2\co QX\rightarrow QD_2X$ gives a stable map
$\Sigma^{\infty}QX\rightarrow \Sigma^{\infty}D_2X$ inducing
a map in homology $h^2_*\co H_nQX \rightarrow H_nD_2X$.
This gives the following commutative diagram.

\begin{equation} \setlength{\unitlength}{3cm}
\begin{picture}(4,0.45)(0,0.45)
\put(1,0.9){\makebox(0,0){$\pi^S_nMO(k)$}}
\put(2,0.9){\makebox(0,0){$\pi^S_nD_2MO(k)$}}
\put(3,0.9){\makebox(0,0){$\pi^S_nMO(2k)$}}
\put(1.35,0.9){\vector(1,0){0.2}}
\put(2.4,0.9){\vector(1,0){0.2}}
\put(1,0){\makebox(0,0){$H_nQMO(k)$}}
\put(2,0){\makebox(0,0){$H_nD_2MO(k)$}}
\put(3,0){\makebox(0,0){$H_nMO(2k)$}}
\put(1.4,0){\vector(1,0){0.2}}
\put(2.4,0){\vector(1,0){0.2}}
\put(1,0.45){\makebox(0,0){$\pi_nQMO(k)$}}
\put(1,0.75){\vector(0,-1){0.2}}
\put(1,0.35){\vector(0,-1){0.2}}
\put(2,0.75){\vector(0,-1){0.6}}
\put(3,0.75){\vector(0,-1){0.6}}
\put(0.95,0.65){\makebox(0,0)[r]{$\cong$}}
\put(1.05,0.25){\makebox(0,0)[l]{$h$}}
\put(2.05,0.45){\makebox(0,0)[l]{$h^S$}}
\put(3.05,0.45){\makebox(0,0)[l]{$h^S$}}
\put(1.45,0.95){\makebox(0,0)[b]{$h^2_*$}}
\put(2.5,0.95){\makebox(0,0)[b]{$\xi_*$}}
\put(1.5,0.05){\makebox(0,0)[b]{$h^2_*$}}
\put(2.5,0.05){\makebox(0,0)[b]{$\xi_*$}}
\end{picture}
\end{equation}
\vspace{1cm}

In this diagram the second and third vertical maps are
stable Hurewicz homormorphisms defined using the fact the
Hurewicz homomorphisms commute with suspension.  The first
square commutes by the definition of the stable Hurewicz
map and by naturality, and the second square commutes by
naturality.

Notice that the normal Stiefel--Whitney numbers (and so the
bordism class) of the double point
self-intersection manifold $\Delta_2(f)$ of an immersion
$f\co M^{n-k}\looparrowright\R^n$ corresponding to
$\alpha\in\pi_n^S MO(k)$ are determined by
(and determine) the Hurewicz image $h^S(\beta)$ of the
element $\beta=\xi_*h^2_*(\alpha)\in\pi_n^S MO(2k)$ corresponding to the
immersion $\theta_2(f)$.  To be more explicit in the case
under consideration we recall the structure of $H_*MO(k)$.

Let $e_i\in H_iBO(1)\cong \Z/2$ be the non-zero element (for
$i\geq 0$).
For each sequence $I=(i_1,i_2,\ldots,i_k)$ of non-negative
integers we define
\[e_I=e_{i_1}e_{i_2}\ldots e_{i_k}=(\mu_k)_*(e_{i_1}\otimes
e_{i_2}\otimes\ldots e_{i_k})\in H_*BO(k)\]
where $\mu_k\co BO(1)^k\rightarrow BO(k)$ is the map
which classifies the product of the universal line bundles.
The dimension of $e_I$ is $|I|=i_1+i_2+\ldots+i_k$.

From the definition of $\mu_k$,
$e_{i_1}e_{i_2}\ldots
e_{i_k}=e_{i_{\sigma(1)}}e_{i_{\sigma(2)}}\ldots
e_{i_{\sigma(k)}}$ for each $\sigma\in\Sigma_k$, where
$\Sigma_k$ is the permutation group on $k$ elements.  Thus each
such element can be written as $e_{i_1}e_{i_2}\ldots
e_{i_k}$ where $i_1\leq i_2\leq\ldots\leq i_k$ and
it follows by a counting
argument that
\[\{\,e_{i_1}e_{i_2}\ldots e_{i_k}\mid 0\leq i_1\leq i_2\leq
\ldots \leq i_k\,\}\]
is a basis for $H_*BO(k)$ (see \cite{koc96}
Proposition~2.4.3).

The sphere bundle of the universal $O(k)$--bundle $\gamma^k$
is given up to homotopy by the inclusion $BO(k-1)\rightarrow
BO(k)$ and so the Thom complex $MO(k)$ is homotopy
equivalent to the quotient space $BO(k)/BO(k-1)$.  It follows
that
\[\{\,e_{i_1}e_{i_2}\ldots e_{i_k}\mid 1\leq i_1\leq i_2\leq
\ldots \leq i_k\,\}\]
is a basis for $\tilde{H}_*MO(k)$.

By Diagrams~(1) and (2),
the double point self-intersection surface of an
immersion $M^{k+2}\looparrowright \R^{2k+2}$ may be
identified up to bordism by using
the stable Hurewicz homomorphism \[h^S\co \pi_{2k+2}^S
MO(2k) \rightarrow H_{2k+2}MO(2k).\]
From the above, $H_{2k+2}MO(2k)$ has a basis
$\{\,e_1^{2k-1}e_3,\; e_1^{2k-2}e_2^2\,\}$.  The element
$\beta\in\pi_{2k+2}^S MO(2k)$ corresponds to an immersion of
a non-boundary (ie,\ a surface of odd Euler characteristic,
bordant to the real projective plane) if
and only if $h^S(\beta)=e_1^{2k-1}e_3$, the only non-zero
stably spherical element.  This corresponds to the fact that
a surface $L$ is a non-boundary if and only if the
normal Stiefel--Whitney number $\overline{w}_1^2[L]=1$ (see
\cite{asa99} Proposition~3.4).

For $k>1$ we are in the stable range,
$\pi_{2k+2}^SMO(2k)\cong\Z/2$, and $h^S$ is a monomorphism:
any two immersions of bordant manifolds are bordant.  On
the other hand, for $k=1$ the group $\Im(2,4)\cong
\pi_4^SMO(2)$ is infinite and $h^S$ is not a monomorphism.
The bordism class of an immersion $L^2\looparrowright
\R^4$ is not determined by $L$.

This discussion can be summed up in the following theorem
which follows essentially from Diagrams~(1) and (2).

\begin{theorem}\label{2.1} Suppose that
$f\co M^{k+2}\looparrowright\R^{2k+2}$ is a
self-transverse immersion corresponding to
$\alpha\in\pi_{2k+2}^S MO(k)$.  Then the double point
self-intersection surface $\Delta_2(f)$ has odd Euler
characteristic and so is a non-boundary if and only if
$$\xi_*h^2_*h(\alpha)=e_1^{2k-1}e_3\in H_{2k+2}MO(2k).\eqno{\qed}$$
\end{theorem}

The map $h^2_*\co H_nQMO(k)\rightarrow H_nD_2MO(k)$
in Diagram~(2) is very easy
to describe in terms of the description of $H_*QX$ as a
Pontrjagin ring
provided by Dyer and Lashof (see \cite{dye62} or
\cite{may76}).  They make use of the
Kudo--Araki operations $Q^i\co H_mQX \rightarrow
H_{m+i}QX$.  These are trivial for $i<m$ and equal to the
Pontrjagin square for $i=m$.  If $I$ denotes the sequence
$(i_1,i_2,\ldots,i_r)$ then we write
$Q^Ix=Q^{i_1}Q^{i_2}\ldots Q^{i_r}x$.  The sequence $I$ is {\em
admissible} if $i_j\leq i_{j+1}$ for $1\leq j<r$ and its
{\em excess} is given by $e(I)=i_1-i_2-\ldots -i_r$.
With this notation we can give the description of $H_*QX$ as
a polynomial algebra:
if $\{\,x_{\lambda}\mid\lambda\in\Lambda\,\}$
is a homogeneous basis for $\tilde{H}_*X\subseteq H_*QX$
where $X$ is a path-connected space then
\[H_*QX=\Z/2[\,Q^Ix_{\lambda}\mid \lambda\in\Lambda,\;I\mbox{ admissible of
excess }e(I)>\dim x_{\lambda}\,].\]
Thus a basis for $H_*QX$ is provided by the monomials in the
polynomial generators.

We may define a height function $\mbox{ht}$ on the
monomial generators of $H_*QX$ by ht$(x_{\lambda})=1$,
ht$(Q^iu)=2\mbox{ht}(u)$ and
ht$(u\cdot v)=$ht$(u)+$ht$(v)$ (where $u\cdot v$ represents
the Pontrjagin product).  The following is a special case of
Lemma~2.3 in \cite{asa99}.

\begin{lemma}\label{2.2} The homomorphism $h^2_*\co
\tilde{H}_*QX \rightarrow
\tilde{H}_*D_2X$ is given by projection onto the monomial
generators of height 2.  The kernel is spanned
by the set of monomials of height other than $2$.
\qed
\end{lemma}

\begin{corollary}\label{2.3} A basis for $H_{2k+2}D_2MO(k)$ is given by
the following set:
$$\{\,e_1^k\cdot e_1^{k-1}e_3,\;e_1^k\cdot e_1^{k-2}e_2^2,\;
e_1^{k-1}e_2\cdot e_1^{k-1}e_2,\;Q^{k+2}e_1^k\,\}.\eqno{\qed}$$
\end{corollary}

The map $\xi_*\co H_nD_2MO(k)\rightarrow H_nMO(2k)$ of
Diagram~(2) is also determined in \cite{asa99}.  In the case
of $n=2k+2$, Theorem~3.1 of \cite{asa99} gives the following.

\begin{lemma}\label{2.4} The homomorphism $\xi_*\co H_{2k+2}D_2MO(k)
\rightarrow H_{2k+2}MO(2k)$ is determined by the following
values:
\begin{enumerate}
\item[] $\xi_*(e_1^k\cdot e_1^{k-1}e_3)=e_1^{2k-1}e_3$;
\item[] $\xi_*(e_1^k\cdot e_1^{k-2}e_2^2)
= \xi_*(e_1^{k-1}e_2\cdot e_1^{k-1}e_2) = e_1^{2k-2}e_2^2$;
\item[] $\xi_*(Q^{k+2}e_1^k)=\left\{\begin{array}{ll}
0 & \mbox{for }k\equiv 0\bmod 4,\\
e_1^{2k-1}e_3 & \mbox{for }k \equiv 1\bmod 4,\\
e_1^{2k-2}e_2^2 & \mbox{for }k\equiv 2\bmod 4,\\
e_1^{2k-1}e_3+e_1^{2k-2}e_2^2 & \mbox{for }k\equiv 3\bmod 4.
\end{array}
\right.$
\end{enumerate}
\end{lemma}

\begin{proof} The first three results are immediate from
\cite{asa99} Theorem~3.1.

For the fourth we apply the formula given in the theorem:
\begin{eqnarray*}
\xi_*(Q^{k+2}e_1^k) & = &
\sum_{m_1+m_2+\ldots+m_k=2}\;\prod_{j=1}^k
{{m_j-1} \choose 0} e_1e_{m_j+1}\quad (m_j\geq 0)\\
& = & \sum_{m_1+m_2+\ldots+m_k=2}
e_1e_{m_1+1}e_1e_{m_2+1}\dots e_1e_{m_k+1}\quad (m_j\geq 0)\\
& = & {k \choose 1}e_1^{2k-1}e_3
+ {k \choose 2}e_1^{2k-2}e_2^2
\end{eqnarray*}
which gives the required result.
\end{proof}

To prove Theorem~\ref{1.1} using Theorem~\ref{2.1} we
determine the image of the spherical classes in
$H_{2k+2}QMO(k)$, ie,\ the classes in the
image of $h\co \pi_{2k+2}QMO(k)\rightarrow H_{2k+2}QMO(k)$,
under $\xi_*\circ h^2_*$ (see Diagram~(2)).
A complete description of these spherical classes is not
necessary for it is sufficient to observe the following
well-known properties of spherical classes (which are
immediate from $H_*S^n$ by naturality).

\begin{lemma}\label{2.5}
{\rm(a)}\qua If an homology class $u\in H_nX$ is spherical then it is
primitive with respect to the cup coproduct,  ie,\
\[\psi(u) = u\otimes 1 + 1 \otimes u,\]
where $\psi\co H_nX
\rightarrow H_n(X\times X)\cong \sum_iH_iX\otimes H_{n-i}X$
is the map induced by the
diagonal map.

{\rm(b)}\qua If an homology class $u\in H_nX$ is spherical (or stably
spherical, ie,\ in the image of $h^S\co \pi^S_nX\rightarrow
H_nX$) then it is annihilated by the reduced Steenrod
algebra, ie,\
\[\mbox{Sq}^i_*(u)=0\]
for all $i>0$, where $\mbox{Sq}^i_*\co H_nX\rightarrow
H_{n-i}X$ is the vector space dual of the usual Steenrod
square cohomology operation $\mbox{Sq}^i\co
H^{n-i}X\rightarrow H^nX$.
\qed
\end{lemma}

To apply the first of these observations we determine the
image of the coproduct primitive submodule of
$H_{2k+2}QMO(k)$ in $H_{2k+2}D_2MO(k)$.

\begin{lemma}\label{2.6}
Suppose that $k>2$.  Then  a
basis for the coproduct primitive
submodule $PH_{2k+2}QMO(k)$ is given by the following set
of elements:
\[\{\,e_1e_{i_2}\ldots e_{i_k}\mid 1\leq
i_2\leq \ldots \leq i_k\,\}
\cup
\{\, e_2^{k-2}e_3^2 + e_1^k\cdot e_1^{k-2}e_2^2 ,\;
e_1^{k-1}e_2\cdot e_1^{k-1}e_2,\;
Q^{k+2}e_1^k\,\}.\]

For $k=1$, a basis for $PH_4QMO(1)$ is given by
\[\{\, Q^{3}e_1,\; e_1\cdot e_1\cdot e_1 \cdot e_1\,\}.\]
For $k=2$, a basis for $PH_6QMO(2)$ is given by
\[\{\, e_1e_5,\; e_3^2+e_1^2\cdot e_2^2+e_1^2\cdot
e_1^2\cdot e_1^2,\; e_1e_2\cdot e_1e_2,\; Q^4e_1^2\,\}.\]
\end{lemma}
\begin{proof}  Recall that $\psi(e_i)=\sum_je_j\otimes
e_{i-j}$.  This determines $\psi(e_I)$ by naturality.

First of all observe that a basis element $e_{i_1}e_{i_2}\ldots
e_{i_k}$ ($i_1\leq i_2\leq \ldots \leq i_k$) of $H_*MO(k)$ is
primitive if and only if $i_1=1$.  To complete
the proof we evaluate the coproduct on the non-primitive
height 1 basis elements and on all the basis elements of greater
height in this
dimension.

For $k>2$ there are no basis elements of height
greater than 2.  For simplicity we use the reduced coproduct
$\tilde\psi(u)=\psi(u)-u\otimes1-1\otimes u$ so that $u$ is
primitive when $\tilde\psi(u)=0$.  Straightforward
calculations give the following:
\begin{enumerate}
\item[] $\tilde\psi(e_2^{k-2}e_3^2)
= e_1^k\otimes e_1^{k-2}e_2^2 + e_1^{k-2}e_2^2\otimes e_1^k$;
\item[] $\tilde\psi(e_2^{k-1}e_4)
= e_1^k\otimes e_1^{k-1}e_3
+ e_1^{k-1}e_2\otimes e_1^{k-1}e_2 + e_1^{k-1}e_3\otimes e_1^k$;
\item[] $\tilde\psi(e_1^k\cdot e_1^{k-1}e_3)
= e_1^k\otimes e_1^{k-1}e_3 + e_1^{k-1}e_3\otimes e_1^k$;
\item[] $\tilde\psi(e_1^k\cdot e_1^{k-2}e_2^2)
= e_1^k\otimes e_1^{k-2}e_2^2 + e_1^{k-2}e_2^2\otimes e_1^k$;
\item[] $\tilde\psi(e_1^{k-1}e_2\cdot e_1^{k-1}e_2)=0$;
\item[] $\tilde\psi(Q^{k+2}e_1^k)=0$.
\end{enumerate}
The lemma follows immediately from these results.

Similar calculations give the results for $k=1$ and $k=2$.
\end{proof}

This lemma has the following immediate corollary.

\begin{corollary}\label{2.7}
For $k>1$, a basis for the projection of
the coproduct primitive submodule $h^2_*PH_{2k+2}QMO(k)
\subseteq H_{2k+2}D_2MO(k)$ is given by the following set
of elements:
\[\{\,e_1^k\cdot e_1^{k-2}e_2^2,\;
e_1^{k-1}e_2\cdot e_1^{k-1}e_2,\;
Q^{k+2}e_1^k\,\}.\]
For $k=1$ a basis for $h_*^2PH_4QMO(1)$ is given by
$\{\,Q^3e_1\,\}$.
\qed
\end{corollary}

To complete the proof of the main theorem it is convenient
to consider various cases depending on the value of $k$
modulo 4.

\section{Case 1: $k$ even}

\begin{theorem}\label{3.1}
For even $k$, the double point
self-intersection surface of each
self-transverse immersion $f\co M^{k+2} \looparrowright
\R^{2k+2}$ has even Euler characteristic and so is a boundary.
\end{theorem}

To prove this result we evaluate dual
Steenrod operations on the elements given by Corollary~\ref{2.7}.

\begin{lemma}\label{3.2}
For even $k$, we have the following results in $H_*D_2MO(k)$:
\begin{enumerate}
\item[] $\Sq^1_*(e_1^k\cdot e_1^{k-2}e_2^2)
= \Sq^1_*(e_1^{k-1}e_2\cdot e_1^{k-1}e_2) = 0$;
\item[] $\Sq^1_*(Q^{k+2}e_1^k) = Q^{k+1}e_1^k$;
\item[] $\Sq^2_*(e_1^k\cdot e_1^{k-2}e_2^2)
= \Sq^2_*(e_1^{k-1}e_2\cdot e_1^{k-1}e_2) = e_1^k\cdot
e_1^k$.
\end{enumerate}
\end{lemma}

\begin{proof}
These are immediate from the Steenrod squares in
$BO(1)=\R P^{\infty}$ ($\Sq^i_*e_j={{j-i}\choose i}e_{j-i}$)
and the Nishida relations (see \cite{may76}).
\end{proof}

\begin{corollary}\label{3.3}
For even $k$, given $\alpha\in\pi_{2k+2}^SMO(k)$,
\[h^2_*h(\alpha)=\lambda(e_1^k\cdot e_1^{k-2}e_2^2
+ e_1^{k-1}e_2\cdot e_1^{k-1}e_2)\]
for $\lambda\in\Z/2$.\end{corollary}

\begin{proof} By Lemma~\ref{2.5}(a) $h^2_*h(\alpha)$ is a
linear combination of the elements given by
Corollary~\ref{2.7}. But since, by Lemma~\ref{2.5}(b), it is
annihiliated by $\Sq^1_*$ and $\Sq^2_*$ only elements of the
given form can arise. \end{proof}

\begin{proof}[Proof of Theorem~\ref{3.1}]
Suppose that $f$ corresponds to $\alpha\in
\pi_{2k+2}QMO(k)$.  Then
\[\xi_*h^2_*h(\alpha) = \lambda\xi_*(e_1^k\cdot
e_1^{k-2}e_2^2 + e_1^{k-1}e_2\cdot e_1^{k-1}e_2)
= \lambda(e_1^{2k-2}e_2^2+e_1^{2k-2}e_2^2)= 0\]
by Lemma~\ref{2.4}.

The result follows by Theorem~\ref{2.1}.
\end{proof}

We can still ask whether there exists an immersion
$f\co M^{k+2}\looparrowright\R^{2k+2}$ corresponding to
$\alpha\in\pi_{2k+2}^SMO(k)$ with $h^2_*h(\alpha)=
e_1^k\cdot e_1^{k-2}e_2^2+e_1^{k-1}e_2\cdot e_1^{k-1}e_2$.
To answer this we first observe the following result.

\begin{proposition}\label{3.4}
For even $k$, given an immersion $f\co
M^{k+2}\looparrowright\R^{2k+2}$ corresponding to
$\alpha\in\pi_{2k+2}^SMO(k)$, the Hurewicz image
$h(\alpha)\in H_{2k+2}QMO(k)$ is determined by the bordism
class of $M$.
\end{proposition}

\begin{proof}
Suppose that $k>2$. Notice that the height 1 part of $h(\alpha)$ is given by
$h^S(\alpha)\in H_{2k+2}MO(k) \subseteq H_{2k+2}QMO(k)$ (see
Diagram~(1) in \cite{asa99}) and so
\[h(\alpha)=h^S(\alpha)+h^2_*h(\alpha)\]
using the obvious inclusion $H_*D_2MO(k)\subseteq
H_*QMO(k)$ as the elements of height 2.
Furthermore, $h^S(\alpha)$ is determined by the normal
Stiefel--Whitney numbers of $M$ and so by the bordism class
of $M$ (see \cite{asa99} Lemma~2.2).

Suppose that $f_1\co M^{k+2}_1\looparrowright\R^{2k+2}$
and $f_2\co M^{k+2}_2\looparrowright\R^{2k+2}$ are two
immersions of bordant manifolds $M_1$ and $M_2$
corresponding to $\alpha_1$, $\alpha_2\in\pi_{2k+2}^SMO(k)$
respectively.  Since $M_1$ and $M_2$ are bordant manifolds
$h^S(\alpha_1)=h^S(\alpha_2)$.  It follows that
\begin{eqnarray*}
h(\alpha_1)-h(\alpha_2) & = & h^2_*(\alpha_1) -
h^2_*(\alpha_2) \\
& = & h^2_*h(\alpha_1-\alpha_2)\\
& = & \lambda(e_1^k\cdot e_1^{k-2}e_2^2 + e_1^{k-1}e_2\cdot
e_1^{k-1}e_2)
\end{eqnarray*}
by Corollary~\ref{3.3}. However $e_1^k\cdot e_1^{k-2}e_2^2 +
e_1^{k-1}e_2\cdot e_1^{k-1}e_2\in H_{2k+2}QMO(k)$ is not
primitive (by Lemma~\ref{2.6}) and so not spherical.

Hence $h(\alpha_1)-h(\alpha_2)=0$ and so
$h(\alpha_1)=h(\alpha_2)$ as required.

For $k=2$ the above proof has to be modified to take account
of height 3 elements.  This leads in this case to
\[h(\alpha_1)-h(\alpha_2)=\lambda(e_1^2\cdot
e_2^2+e_1e_2\cdot e_1e_2 + e_1^2\cdot e_1^2\cdot e_1^2)\]
which gives the same result.
\end{proof}

From now on we write $\alpha(m)$ for the number of digits $1$ in
the dyadic expression for the positive integer $m$.

\begin{proposition}\label{3.5}
For even $k>2$, given an immersion $f\co
M^{k+2}\looparrowright \R^{2k+2}$ corresponding to
$\alpha\in\pi_{2k+2}^SMO(k)$, the Hurewicz image $h(\alpha)\in
H_{2k+2}QMO(k)$ is given as follows:
\[ h(\alpha) = \left\{\begin{array}{l}
h^S(\alpha),\quad\mbox{if the normal Stiefel--Whitney number }
\overline{w}_2\overline{w}_k[M]=0,\\
h^S(\alpha) + e_1^k\cdot e_1^{k-2}e_2^2+e_1^{k-1}e_2\cdot
e_1^{k-1}e_2,\quad\mbox{if\/ }
\overline{w}_2\overline{w}_k[M]=1.\end{array}\right.\]
The second case arises if and only if $\alpha(k+2)\leq 2$.
In this case, writing $k+2=2^r+2^s$ where $r\geq s\geq 1$, we
can take $M$ to be $\R P^{2^r}\times \R P^{2^s}$.

For $k=2$, given $f\co M^4\looparrowright \R^6$
corresponding to $\alpha\in\pi_6^SMO(2)$, the Hurewicz image
$h(\alpha)\in H_6QMO(2)$ is given as follows:
\[ h(\alpha) = \left\{\begin{array}{l}
0,\quad\mbox{if\/ $\overline{w}_2^2[M]=0$,}\\
e_3^2+e_1^2\cdot e_2^2+e_1e_2\cdot e_1e_2 + e_1^2\cdot
e_1^2\cdot e_1^2,\quad\mbox{if\/ $\overline{w}_2^2[M]=1$.}
\end{array}\right.\]
\end{proposition}

\begin{proof}  Suppose that $k>2$ and that $f\co M^{k+2}\looparrowright
\R^{2k+2}$ is an immersion corresponding to
$\alpha\in\pi_{2k+2}^SMO(k)$.

If $\lambda=0$ in Corollary~\ref{3.3} then
$h(\alpha)=h^S(\alpha)$.  By Lemma~\ref{2.6}, $\lambda=1$ if
and only if $e_2^{k-2}e_3^2$ has coefficient $1$ when
$h^S(\alpha)$ is written in terms of the basis $\{e_I\}$.
This occurs if and only if the Stiefel--Whitney number
$\overline{w}_2\overline{w}_k[M]=1$. For, by \cite{asa99}
Lemma~2.2, $\overline{w}_2\overline{w}_k[M]$ is given by the
Kronecker product $\langle w_2w_k^2,h^S(\alpha)\rangle$ in
$MO(k)$.  By naturality, this product can be evaluated in
$BO(1)^k$ using $\mu_k$ and the
fact that
$\mu_k^*w_i=\sigma_i(x_1,x_2,\ldots,x_k)\in\Z/2[x_1,x_2,\ldots
x_k]\equiv H^*BO(1)^k$, the $i$th  elementary symmetric
polynomial.  In this case $\mu_k^*(w_2w_k^2)=\sum
x_1^3 x_2^3 x_2^2 \dots x_k^2$ and so the Kronecker product $\langle
w_2w_k^2,h^S(\alpha)\rangle$ is given by the coefficient of
$e_2^{k-2}e_3^2$.  Hence $\lambda=1$ if and only if
$\overline{w}_2\overline{w}_k[M]=1$.

In the case $\alpha(k+2)>2$, by a
theorem of R\,L\,W~Brown (\cite{bro71} Theorem~5.1), there exists an
embedding $f_1\co M_1\hookrightarrow \R^{2k+2}$ of a
manifold $M_1$ bordant to $M$.  Suppose that $f_1$
corresponds to $\alpha_1\in\pi_{2k+2}^SMO(k)$.  Then
$h^2_*h(\alpha_1)=0$ since$f_1$ is an embedding with no
double points.  Hence, by Proposition~3.4,
$h^2_*h(\alpha)=0$ and so $h(\alpha)=h^S(\alpha)$.

For $\alpha(k+2)\leq 2$ put $k+2=2^r+2^s$ where $r\geq s\geq
1$.  Then by Whitney's immersion theorem there are
immersions $\R P^{2^r}\looparrowright \R^{2^{r+1}-1}$ and
$\R P^{2^s}\looparrowright \R^{2^{s+1}-1}$ and the product
of these gives an immersion $f\co \R P^{2^r}\times \R
P^{2^s}\looparrowright \R^{2^{r+1}+2^{s+1}-2}=\R^{2k+2}$.
Finally a standard verification shows that the normal
Stiefel--Whitney number
$\overline{w}_2\overline{w}_{2^r+2^s-2}[\R P^{2^r}\times\R
P^{2^s}]=1$.

The argument for $k=2$ is almost identical.  The manifold
$M^4$ either is a boundary or is bordant to $\R P^2\times
\R P^2$ depending on the value of the normal Stiefel--Whitney
number $\overline{w}_2^2[M]$.  The presence of the height~3 term
$e_1^2\cdot e_1^2\cdot e_1^2$ in the Hurewicz image shows
that in this case any immersion of a non-boundary has an odd
number of triple points.
\end{proof}

\section{Case 2: $k\equiv 1\bmod 4$}

\begin{theorem}\label{4.1}
For $k\equiv 1\bmod 4$, given any manifold
$M^{k+2}$ there exists an immersion $M^{k+2}\looparrowright
\R^{2k+2}$ with double point self-intersection surface of odd Euler
characteristic (and so a non-boundary)
and there exists another immersion with
double point self-intersection surface of even Euler
characteristic (and so a boundary).
\end{theorem}

\begin{proof} For $k\equiv 1\bmod 4$, $\alpha(k+2)\geq 2$
and so there exists an immersion $f\co M^{k+2}\looparrowright
\R^{2k+2}$ by R~Cohen's immersion theorem (\cite{coh85}).
Furthermore, for such $k$ there is an immersion
$S^{k+2}\looparrowright\R^{2k+2}$ with double point manifold
of odd Euler characteristic (\cite{ecc97} Theorem~1.2).
Taking the connected sum of this immersion and $f$ gives an
immersion of $M$ with a double point self-intersection surface
with Euler characteristic of opposite parity to that of the
double point self-intersection surface of $f$.  Hence both
parities can arise.
\end{proof}

Again it is natural to ask about the Hurewicz image
$h(\alpha)\in H_{2k+2}MO(k)$ of the elements
$\alpha\in\pi_{2k+2}^SMO(k)$ corresponding to these
immersions.  The situation is the following.

\begin{proposition}\label{4.2}
For $k\equiv 1\bmod 4$ such that $k>1$, given an immersion $f\co
M^{k+2}\looparrowright \R^{2k+2}$ corresponding to
$\alpha\in\pi_{2k+2}^SMO(k)$, then the Hurewicz image $h(\alpha)\in
H_{2k+2}QMO(k)$ is given as follows:
\[h(\alpha)=
\left\{\begin{array}{l}
h^S(\alpha),\quad\mbox{if $\Delta_2(f)$ has even Euler
characteristic,}\\
h^S(\alpha)+Q^{k+2}e_1^k,\quad\mbox{if $\Delta_2(f)$ has odd Euler
characteristic.}\end{array}\right.\]

For $k=1$, given an immersion $f\co
M^3\looparrowright \R^4$ corresponding to
$\alpha\in\pi_4^SMO(1)$, then the Hurewicz image $h(\alpha)\in
H_4QMO(1)$ is given as follows:
\[h(\alpha)=
\left\{\begin{array}{l}
0,\quad\mbox{if $\Delta_2(f)$ has even Euler
characteristic,}\\
Q^{3}e_1+e_1\cdot e_1\cdot e_1\cdot e_1,
\quad\mbox{if $\Delta_2(f)$ has odd Euler
characteristic.}\end{array}\right.\]
\end{proposition}

\begin{proof}
Evaluating the dual Steenrod operations on
the elements of Corollary~\ref{2.7} we obtain
\begin{enumerate}
\item[] $\Sq^1_*(e_1^k\cdot e_1^{k-2}e_2^2)
=\Sq^1_*(e_1^{k-1}e_2\cdot e_1^{k-1}e_2)
=\Sq^1_*(Q^{k+2}e_1^k) = 0$;
\item[] $\Sq^2_*(e_1^k\cdot e_1^{k-2}e_2^2)
=\Sq^2_*(e_1^{k-1}e_2\cdot e_1^{k-1}e_2)
=e_1^k\cdot e_1^k$;
\item[] $\Sq^2_*(Q^{k+2}e_1^k)=0$.
\end{enumerate}
Hence $h^2_*h(\alpha)$ is in the submodule spanned by the
set
\[\{\,e_1^k\cdot e_1^{k-2}e_2^2+e_1^{k-1}e_2\cdot
e_1^{k-1}e_2,\;Q^{k+2}e_1^k\,\}.\]

For $k\equiv 1\bmod
4$ and $k>1$ it follows that $\alpha(k) > 2$ and so, by
R\,L\,W~Brown's embedding theorem (\cite{bro71} Theorem~5.1), $M$ is
bordant to a manifold $M_1$ which has an embedding
$f_1\co M_1\hookrightarrow\R^{2k+2}$.  Let $f_1$
correspond to $\alpha_1\in\pi_{2k+2}^SMO(k)$.
Then $h^2_*h(\alpha_1)=0$ since $f_1$ is an
embedding.  Furthermore, $h^S(\alpha)=h^S(\alpha_1)$ since
these are determined by the bordism class.  It
follows, by Lemma~\ref{2.6},
that $h^2_*h(\alpha)=h(\alpha-\alpha_1)=\lambda
Q^{k+2}e_1^k$ since this is necessarily a primitive class in
$H_{2k+2}QMO(k)$.  Hence $h(\alpha)=h^S(\alpha)+\lambda
Q^{2k+2}e_1^k$ where, by Theorem~\ref{2.1} and Lemma~\ref{2.4},
$\lambda$ gives
the parity of the Euler characteristic of the double point
self-intersection surface of the immersion.

For $k=1$, given any immersion $f\co M^3\looparrowright\R^4$
the manifold $M$ is necessarily a boundary since this is
true of all 3--manifolds and so, if
$\alpha\in\pi_4^SMO(1)\cong\Z/2$ is the corresponding
element, $h^S(\alpha)=0$.  The formula for $h(\alpha)$ then
follows from the value for the Hurewicz image coming from
\cite{ecc80} Proposition~3.4.  The presence of the height~4
term $e_1\cdot e_1\cdot e_1\cdot e_1$ when $\Delta_2(f)$ has
odd Euler characteristic indicates that such immersions have
an odd number of quadruple points.
\end{proof}

\section{Case 3: $k\equiv 3 \bmod4$}

\begin{theorem}\label{5.1}
For $k\equiv 3 \bmod 4$, there exists a manifold $M^{k+2}$
with a self-transverse
immersion $M^{k+2}\looparrowright \R^{2k+2}$
for which the Euler characteristic of the
double point self-intersection surface is odd (so that it is
a non-boundary)  if and only if
$k+1$ is a power of $2$.  In this case, the parity of the
Euler characteristic of the double point self-intersection
surface of a self-transverse immersion $M^{k+2}\looparrowright \R^{2k+2}$
is given by the normal Stiefel--Whitney number
$\overline{w}_2\overline{w}_k[M]$.
\end{theorem}

To prove this we begin as in Case~1.

\begin{lemma}\label{5.2}
For $k\equiv 3 \bmod 4$, we have the following results in
$H_*D_2MO(k)$:
\begin{enumerate}
\item[] $\Sq^1_*(e_1^k\cdot e_1^{k-2}e_2^2)
=\Sq^1_*(e_1^{k-1}e_2\cdot e_1^{k-1}e_2)
=\Sq^1_*(Q^{k+2}e_1^k)
=0$;
\item[] $\Sq^2_*(e_1^k\cdot e_1^{k-2}e_2^2)
=\Sq^2_*(e_1^{k-1}e_2\cdot e_1^{k-1}e_2)
=\Sq^2_*(Q^{k+2}e_1^k)
=e_1^k\cdot e_1^k$.\qed
\end{enumerate}
\end{lemma}

\begin{corollary}\label{5.3}
For $k\equiv 3 \bmod 4$, given
$\alpha\in\pi_{2k+2}^SMO(k)$, $h^2_*h(\alpha)$ lies in the
submodule of $H_{2k+2}D_2MO(k)$ spanned by the following
set:
$$\{\,
e_1^k\cdot e_1^{k-2}e_2^2 + e_1^{k-1}e_2\cdot
e_1^{k-1}e_2,\;
e_1^k\cdot e_1^{k-2}e_2^2 + Q^{k+2}e_1^k,\,\}.\eqno{\qed}$$
\end{corollary}

To prove Theorem~\ref{5.1} we
determine which of the elements in the submodule of
Corollary~\ref{5.3}
can arise as a value of $h^2_*h$.  First of all we observe that, for
$\alpha\in\pi_{2k+2}^SMO(k)$ corresponding to a
self-transverse immersion $f\co
M^{k+2}\looparrowright\R^{2k+2}$, the value of
$h^2_*h(\alpha)$ and so the bordism class of the double point
self-intersection surface is determined by the bordism class
of $M$.
This is the analogue of Proposition~\ref{3.4}.

\begin{proposition}\label{5.4}
For $k \equiv 3 \bmod 4$, given an immersion $f\co
M^{k+2}\looparrowright\R^{2k+2}$ corresponding to
$\alpha\in\pi_{2k+2}^SMO(k)$, the Hurewicz image
$h(\alpha)\in H_{2k+2}QMO(k)$ is determined by the bordism
class of $M$.
\end{proposition}

Before proving this we use it to complete the proof of
Theorem~\ref{5.1} for $k+1$ not a power of $2$.

\begin{corollary}\label{5.4a}
For $k \equiv 3 \bmod 4$, given a self-transverse immersion $f\co
M^{k+2}\looparrowright\R^{2k+2}$, the bordism class of
$\Delta_2(f)$, the double point self-intersection surface,
is determined by the bordism class of $M$.
\end{corollary}

\begin{proof}
This is immediate from the preceding proposition since, by
Theorem~\ref{2.1}, the
bordism class of $\Delta_2(f)$ is determined by
$h^2_*h(\alpha)$.
\end{proof}

Theorem~\ref{5.1} for $k+1$ not a power of $2$ is an immediate
consequence of this.

\begin{corollary}\label{5.4b}
For $k \equiv 3 \bmod 4$ such that $\alpha(k+2)>2$, the
double point self-intersection surface of each
self-transverse immersion $f\co
M^{k+2}\looparrowright\R^{2k+2}$ has even Euler
characteristic and so is a boundary.
\end{corollary}

\begin{proof}
By R\,L\,W~Brown's embedding theorem (\cite{bro71}
Theorem~5.1), $M$ is
bordant to a manifold $M_1$ which has an embedding
$f_1\co M_1 \hookrightarrow \R^{2k+2}$.  The double point
self-intersection surface of $f_1$ is empty.  Hence, by
Corollary~\ref{5.4a}, the double point self-intersection
surface of $f$ is a boundary.
\end{proof}

We shall return to Theorem~\ref{5.1} for $k$ such that $\alpha(k+2)=2$,
ie,\ $k+1$ is a power of 2,
after proving Proposition~\ref{5.4}.

\begin{proof}[Proof of Proposition~\ref{5.4}]
Suppose that $f_1\co M^{k+2}_1\looparrowright\R^{2k+2}$
and $f_2\co M^{k+2}_2\looparrowright\R^{2k+2}$ are two
immersions of bordant manifolds
corresponding to $\alpha_1$, $\alpha_2\in\pi_{2k+2}^SMO(k)$
respectively.  As in the proof of Proposition~\ref{3.4},
since $M_1$ and $M_2$ are bordant manifolds
$h^S(\alpha_1)=h^S(\alpha_2)$.  It follows that
\[h(\alpha_1-\alpha_2) = h(\alpha_1)-h(\alpha_2) =
h^2_*h(\alpha_1)-h^2_*h(\alpha_2) =
h^2_*h(\alpha_1-\alpha_2)\]
and so lies in the submodule of $H_{2k+2}QMO(k)$ spanned
by the set given by Corollary~\ref{5.3}.

We consider the non-zero elements of this submodule in turn.
First of all, by Lemma~\ref{2.6}, neither
$e_1^k\cdot e_1^{k-2}e_2^2 + e_1^{k-1}e_2\cdot e_1^{k-1}e_2$
nor
$e_1^k\cdot e_1^{k-2}e_2^2 + Q^{k+2}e_1^k$
is primitive and so neither of these elements can be
spherical.
This leaves the element $e_1^{k-1}e_2\cdot e_1^{k-1}e_2 +
Q^{k+2}e_1^k$, which is primitive.

\begin{lemma}\label{5.5}
The element $e_1^{k-1}e_2\cdot e_1^{k-1}e_2 +
Q^{k+2}e_1^k \in H_{2k+2}QMO(k)$ is not spherical.
\end{lemma}

Accepting this for the moment, it follows that
$h(\alpha_1)-h(\alpha_2)
=h(\alpha_1-\alpha_2)=0$ and so
$h(\alpha_1)=h(\alpha_2)$ as required.
\end{proof}

The proof of Lemma~\ref{5.5} needs an additional idea.
Certain elements
$\alpha\in\pi_nQX\cong \pi_n^SX$ are detected
by their Hurewicz image in $H_nQX$.  Elements may also be
detected by a
cohomology operation in $H^*C_{\alpha}$ where the stable
space $C_{\alpha}$ is the mapping cone of the stable map
$\alpha\co S^n\rightarrow X$.  It is interesting to
consider the relationship between these two approaches.  One
of the most basic facts is the following.

\begin{proposition}\label{5.6}
Given $\alpha\in \pi_{2m}QX\cong \pi_{2m}^SX$, if
$h(\alpha)=u^2\in H_{2m}QX$ where $u\in H_mX$, then
$\Sq^{m+1}\overline{u}\not= 0\in H^{2m+1}C_{\alpha}$, where
$\overline{u}\in H^mC_{\alpha}$ is any class such that the
Kronecker product $\langle i^*\overline{u},u\rangle=1$
writing $i\co X\rightarrow C_{\alpha}$ for the natural inclusion
map.
\end{proposition}

Notice that
$i^*\Sq^{m+1}\overline{u}=\Sq^{m+1}i^*\overline{u}=0\in H^{2m+1}X$ for
dimensional reasons and so $\Sq^{m+1}\overline{u}$ lies in
$p^*H^{2m+1}S^{2m+1}\cong\Z/2$ where $p\co
C_{\alpha}\rightarrow S^{2m+1}$ is the natural projection
map.  When $\Sq^{m+1}\overline{u}\not=0$
the stable space $C_{\alpha}$ requires at least one
suspension before it can exist as an unstable space.

\begin{proof}[Proof of Proposition~\ref{5.6}]
This proposition can be thought of as describing how
detection by a certain
James--Hopf invariant corresponds to detection by a Steenrod--Hopf
invariant.  The proposition follows from work of Boardman and Steer
(\cite{boa67} Corollary 5.15).  The special case of $X=\R
P^{\infty}$ appears as \cite{ecc81} Lemma~4.2 and the
justification given there extends to the general case.
\end{proof}

\begin{proof}[Proof of Lemma~\ref{5.5}]
Suppose, for contradiction, that there exists an element
$\alpha\in \pi_{2k+2} QMO(k)$ such that
\[h(\alpha) =
e_1^{k-1}e_2\cdot e_1^{k-1}e_2 + Q^{k+2}e_1^k \in
H_{2k+2}QMO(k).\]
Adjointness gives a natural isomorphism $\sigma\co
\pi_nQX \cong \pi_{n+1}Q\Sigma X$.  Since the corresponding homology
suspension $\sigma\co H_nQX\rightarrow H_{n+1}Q\Sigma X$
kills Pontrjagin products and commutes with the Kudo--Araki
operations,
\[h(\sigma^2\alpha) = Q^{k+2}\sigma^2e_1^k =
\sigma^2e_1^k\cdot \sigma^2e_1^k \in H_{2k+4}
Q\Sigma^2MO(k).\]
Hence, by Proposition~\ref{5.6},
\begin{equation}\label{eqn3}
\Sq^{k+3}\overline{u}\not=0\in
H^{2k+5}C_{\sigma^2\alpha}
\end{equation}
where $\overline{u}\in H^{k+2}C_{\sigma^2\alpha}$ is an
element such that $\langle i^*\overline{u},\sigma^2e_1^k\rangle=1$.

We now show that this is impossible.

First of all notice that, since $h^S(\sigma^2\alpha)=0\in
H_{2k+4}\Sigma^2MO(k)$,
$i^*\co H^iC_{\sigma^2\alpha}\cong H^i\Sigma^2MO(k)$ for
$i\not=2k+5$.  This means that away from dimension $2k+5$ we
can do calculations in $H^*C_{\sigma^2\alpha}$ by calculating
in $H^*\Sigma^2MO(2)$.
In particular, $H^{k+2}C_{\sigma^2\alpha}\cong
H^{k+2}\Sigma^2MO(k) \cong\Z/2$ is generated by $\sigma^2
w_k$ and so $\overline{u}=\sigma^2 w_k$. Here $w_k\in
H^kMO(k)\subseteq H^k BO(k)$ is the universal
Stiefel--Whitney class.

Secondly, we can evaluate the action of the Steenrod
squaring operations on $H^*\Sigma^2MO(k)$ using the fact that
the operations commute with suspension and using their action on
$H^*MO(k)$ which is determined by the Wu
formula (see \cite{mil74} Problem~8A).
In evaluating the operations we use the fact that
$k\equiv 3\bmod 4$ which means that we can write $k=4r-1$
for some positive
integer $r$.

Then,
\begin{eqnarray*}
\Sq^{k+3}\overline{u} & = & \Sq^{4r+2}\sigma^2w_{4r-1}\\
& = &
(\Sq^2\Sq^{4r}+\Sq^1\Sq^{4r}\Sq^1)\sigma^2w_{4r-1}\quad
\mbox{by the Adem relations}\\
& = & \Sq^2\sigma^2\Sq^{4r}w_{4r-1} +
\Sq^1\sigma^2\Sq^{4r}\Sq^1w_{4r-1}\\
&& \qquad\qquad
\mbox{since the operations commute with suspension}\\
& = & 0 + \Sq^1\sigma^2w_1^2w_{4r-1}^2\quad
\mbox{by the Wu formula and dimension}\\
& = & \Sq^1\sigma^2\Sq^1w_1w_{4r-1}^2\quad
\mbox{by a further application of the Wu formula}\\
& = & \Sq^1\Sq^1\sigma^2w_1w_{4r-1}^2\quad
\mbox{since $\Sq^1$ commutes with suspension}\\
& = & 0 \quad
\mbox{by the Adem relations.}
\end{eqnarray*}
This contradicts equation~(\ref{eqn3}) and so proves the
lemma, completing the proof of Proposition~\ref{5.4} and so
the proof of Theorem~\ref{5.1} in the case of $\alpha(k+2)>2$.
\end{proof}

Now suppose that $k\equiv 3\bmod 4$ and $\alpha(k+2)=2$.  In
this case $k+1=2^r$ for some integer $r$ such that $r\geq
2$.

To complete the proof of Theorem~\ref{5.1} we consider
the structure of the bordism group $MO_{k+2}$ of
$(k+2)$--dimensional manifolds.  Recall that the unoriented
bordism ring $MO_*$
is a polynomial ring with one generator in each dimension
not of the form $2^s-1$.  We describe a manifold $M^{k+2}$
as decomposable (up to bordism) if it represents a
decomposable element in this ring and so is bordant to a union of
product manifolds $N_1^{m_1}\times N_2^{m_2}$ where $m_1$,
$m_2>0$.  We can assume that $m_1$, $m_2\geq 2$ since there
is no non-boundary in dimension 1.  Any indecomposable
manifold of dimension $k+2$ can be taken as a representative of a polynomial
generator in this dimension and any two such manifolds are
bordant modulo decomposable manifolds.  Notice that, by
R~Cohen's immersion theorem (\cite{coh85}), every manifold in
this dimension immerses in $\R^{2k+2}$.

We first of all deal with the decomposable manifolds.

\begin{proposition}\label{5.9}
For $k\equiv 3\bmod 4$ such that $\alpha(k+2)=2$, every
manifold $M^{k+2}$ which is decomposable in the bordism ring
is bordant to a manifold which embeds in $\R^{2k+2}$.  Hence
the double point self-intersection surface of any immersion
$M^{k+2}\looparrowright\R^{2k+2}$ has even Euler
characteristic and so is a boundary.
\end{proposition}

\begin{proof} To prove that decomposable manifolds of
dimension $k+2$ embed
up to bordism it is sufficient to prove that each product
$N_1^{m_1}\times N_2^{m_2}$, with $m_1$, $m_2\geq 2$ and
$m_1+m_2=k+2$, embeds up to bordism in $\R^{2k+2}$.  To do this we make
use of the following result which follows easily from the
Whitney embedding theorem.

\begin{lemma}[\cite{bro71} Lemma~2.1]\label{5.10}
If the manifold $N_1^{m_1}$ immerses in $\R^s$, the manifold
$N_2^{n_2}$ embeds in $\R^t$ and $s+t\geq 2m_1+1$ then
$N_1\times N_2$ embeds in $\R^{s+t}$.
\qed
\end{lemma}

Since $m_1+m_2=k+2=2^r+1$ and $m_1$, $m_2\geq 2$ it follows
that $\alpha(m_1)\geq 2$ or $\alpha(m_2)\geq 2$.  Suppose
without loss of generality that $\alpha(m_2)\geq 2$.  Then,
by R\,L\,W~Brown's embedding theorem
(\cite{bro71}), $N_2$ is bordant to a manifold which embeds
in $\R^{2m_2-1}$.  In addition,
by the Whitney immersion theorem, $N_1$ immerses in
$\R^{2m_1-1}$.  Hence, by the lemma, $N_1\times N_2$ is
bordant to a manifold which embeds in
$\R^{2m_1+2m_2-2}=\R^{2k+2}$.  The numerical condition of
the lemma is automatically satisfied: $s+t=2k+2\geq 2m_1+1$
since $m_1\leq k$.

The final part of the proposition is immediate from
Corollary~\ref{5.4a} as in the proof of
Corollary~\ref{5.4b}.
\end{proof}

We turn now to the indecomposable manifolds. 

\begin{proposition}\label{5.11}
For $k\equiv 3 \bmod 4$ such that $\alpha(k+2)=2$, suppose
that
$\alpha\in\pi_{2k+2}^SMO(k)$ corresponds to an immersion
$f\co M^{k+2}\looparrowright \R^{2k+2}$ of a manifold $M$
which is indecomposable in the bordism ring.  Then
\[h(\alpha)=h^S(\alpha)+ e_1^k\cdot e_1^{k-2}e_2^2 +
Q^{k+2}e_1^k\]
and so the double point self-intersection surface
$\Delta_2(f)$ has odd Euler characteristic and so is not a
boundary.
\end{proposition}

\begin{proof}
We make use of particular indecomposable manifolds
constructed by A~Dold (\cite{dol56}).  In dimension
$k+2=2^r+1$ this manifold $V^{k+2}$ is formed from the
product $S^1\times \C P^{2^{r-1}}$ by identifying $(u,z)$
with $(-u,\overline{z})$. In Dold's notation this is
$P(1,2^{r-1})$.  He shows that the cohomology
ring of $V$ is given by
\[H^*V = \Z/2[c,d]/(c^2,d^{2^{r-1}+1})\]
where $\dim(c)=1$ and $\dim(d)=2$, and the total tangent
Stiefel--Whitney class of $V$ is given by
\[ w(V) = (1+c)(1+c+d)^{2^{r-1}+1}.\]
This implies that the total normal Steifel--Whitney class of
$V$ is given by
\[ \overline{w}(V) = (1+c)(1+c+d)^{2^{r-1}-1}\]
from which it follows that $\overline{w}_1(V)=0$,
$\overline{w}_2(V)=d$ and
$\overline{w}_k=\overline{w}_{2^r-1}=cd^{2^{r-1}-1}$.

Hence $V$ is an orientable manifold such that the normal
Stiefel--Whitney number $\overline{w}_2\overline{w}_k[V]=1$.

Let $f:V\looparrowright\R^{2k+2}$ be an immersion
corresponding to $\alpha\in\pi_{2k+2}^SMO(k)$.
Using the argument in the proof if Proposition~\ref{3.5} it
follows from the non-vanishing of this Stiefel--Whitney number
that $e_2^{k-2}e_3^2$ has coefficient 1 when
$h^S(\alpha)\in H_{2k+2}MO(k)$ is written in terms of the
basis $\{e_I\}$.  Hence from Lemma~\ref{2.6} and
Corollary~\ref{5.3} it follows that either
\[ h(\alpha) = h^S(\alpha) + e_1^k\cdot e_1^{k-2}e_2^2 +
e_1^{k-1}e_2\cdot e_1^{k-1}e_2 \]
or
\[ h(\alpha) = h^S(\alpha) + e_1^k\cdot e_1^{k-2}e_2^2 +
Q^{k+2}e_1^k.\]

However, since $V$ is orientable, $h(\alpha)$ lies in the
image of $H_{2k+2}QMSO(k)$ in $H_{2k+1}QMO(k)$.  Since
$H_{k+1}MSO(k)=0$ the element $e_1^{k-1}e_2$ does not come
from $H_{2k+1}MSO(k)$ which eliminates the first of the
above possibilities.  This proves the formula in the
proposition for the manifold $V$.  However, since any other
indecomposable manifold $M$ is bordant to $V$ modulo a
decomposable manifold the same formula holds for an
immersion of $M$ by Proposition~\ref{5.4} and
Proposition~\ref{5.9}.

The Euler characteristic of the double point
self-intersection surface is now given by Theorem~\ref{2.1}
since, by Lemma~\ref{2.4}, $\xi_*(e_1^k\cdot e_1^{k-2}e_2^2 +
Q^{k+2}e_1^k)=e_1^{2k-1}e_3$.
\end{proof}

\begin{proof}[Proof of Theorem~\ref{5.1}]
Almost everything is given by Corollary~\ref{5.4b},
Proposition~\ref{5.9} and Proposition~\ref{5.11}.  The final
observation about the Stiefel--Whitney number is clear since,
by Lemma~\ref{2.6}, $e_1^k\cdot e_1^{k-2}e_2^2$ occurs in
$h(\alpha)$ if and only if $e_2^{k-2}e_3^2$ does and
by the argument in the proof of Proposition~\ref{3.5} this
is equivalent to $\overline{w}_2\overline{w}_k[M]=1$.
\end{proof}

\section{Comparison with Sz\H{u}cs's results}

To conclude, we compare the results and methods of this paper
with those of Sz\H{u}cs in \cite{szu93}. His methods are more
geometric than ours and he gives explicit constructions to show that
any surface of even Euler characteristic can occur as the double
point self-intersection manifold of a self-tranverse immersion
$M^{k+2}\looparrowright\R^{2k+2}$ if $k$ is even, and any surface can
arise if $k\equiv 1\bmod 4$.

There are two key steps in Sz\H{u}cs's non-existence proof.

\begin{proposition}[\cite{szu93} Claim~1]\label{6.1}
For $k\not\equiv 1 \bmod 4$, the parity of the Euler
characteristic of the double point self-intersection surface
of a self-transverse immersion
$M^{k+2}\looparrowright\R^{2k+2}$
depends only on the bordism class of $M$.
\end{proposition}

\begin{proposition}[\cite{szu93} Claim~3]\label{6.2}
For $k\not\equiv 1 \bmod 4$, the Euler
characteristic of the double point self-intersection surface
of a self-transverse immersion
$M^{k+2}\looparrowright\R^{2k+2}$ is even, if $M$ is
decomposable in the bordism ring.
\end{proposition}

These are both included in our results. In fact we
have proved a stronger result than Proposition~\ref{6.1} as
follows.

\begin{proposition}\label{6.3}
For $k\not\equiv 1 \bmod 4$,  given an element
$\alpha\in\pi_{2k+2}^SMO(k)$ corresponding to
a self-transverse immersion
$M^{k+2}\looparrowright\R^{2k+2}$,
the Hurewicz image
$h(\alpha)\in H_{2k+2}QMO(k)$
depends only on the bordism class of $M$.
\end{proposition}

This is our Proposition~\ref{3.4} and Proposition~\ref{5.4}.
It implies Proposition~\ref{6.1} by Theorem~\ref{2.1}.

Sz\H{u}cs's technique for proving Proposition~\ref{6.1} is
to extend the discussion to include maps which he calls {\em
prim} maps; these are smooth maps
$M^{k+2}\rightarrow \R^{2k+2}$ which arise as the projection
of an immersion
$M^{k+2}\looparrowright\R^{2k+3}\rightarrow\R^{2k+2}$.
Sz\H{u}cs claims that Proposition~\ref{6.1} remains true for
these maps.  However, this must be false since he describes
(\cite{szu93} Lemma~4)
how to construct a prim map $f\co V^{2^r+1}\rightarrow
\R^{2^{r+1}}$ of an indecomposable manifold with a double
point self-intersection surface of even Euler
characteristic, whereas by our Proposition~\ref{5.11} the double point
self-intersection surface of any immersion of $V$ will have
odd Euler characteristic.

Turning to Proposition~\ref{6.2},
we have proved this explicitly for $k\equiv 3\bmod 4$ as
Proposition~\ref{5.9}.  It is immediate for even $k$ from
Theorem~\ref{3.1}.  It is intriguing that, whereas by our methods the
case of even $k$ is the simplest, by Sz\H{u}cs's
geometric approach the argument in this case is quite delicate.  He
readily reduces the question to the case of $\R
P^{2^r}\times \R P^{2^s}$ (using Proposition~\ref{6.1} and
Brown's embedding theorem as we do) but then dealing with these
manifolds calls for some ingenuity.  It is not surprising
that these manifolds are the difficult ones since these are
the even dimensional manifolds for which $h^2_*h(\alpha)$ is
non-zero (by Proposition~\ref{3.5}).

\hyphenation{Whit-ney}

\end{document}